\documentclass[12pt,dvips]{amsart}
\usepackage{euler, amsfonts, amssymb, latexsym, epsfig,epic}

\usepackage{amscd,amssymb}
\usepackage{verbatim}
\usepackage{pstricks}
\usepackage{pst-node}

\newcommand\lie[1]{{\mathfrak #1}}

\newcommand\diag{{\rm diag}}

\newcommand\iso{\mathop{\cong}}
\newcommand\tensor{{\otimes}}

\newtheorem{Theorem}{Theorem}
\newtheorem{Proposition}{Proposition}
\newtheorem{Lemma}{Lemma}

\newtheorem{Corollary}{Corollary}
\newtheorem*{Corollary*}{Corollary}

\newtheorem*{Theorem*}{Theorem}
\newtheorem*{Remark}{Remark}
\theoremstyle{remark}
\newtheorem{Example}{Example}

\newcommand\onto{\mathop{\twoheadrightarrow}}
\newcommand\into{\operatorname*{\hookrightarrow}}
\newcommand\To{\longrightarrow}
\newcommand\From{\longleftarrow}

\newcommand\Pone{{\mathbb P}^1}
\newcommand\PP{\mathbb P}

\newcommand\complexes{{\mathbb C}}

\newcommand\naturals{{\mathbb N}}

\theoremstyle{plain}

\newcommand\dfn{\bf} 

\newtheorem{thm}{Theorem}

\theoremstyle{definition}

\theoremstyle{remark}
\newtheorem{rem}[thm]{Remark}
\newtheorem{ex}[thm]{Example}

\newcommand{\Z}{\mathbb{Z}}
\newcommand{\C}{\mathbb{C}}

\newcommand{\Ext}{\mathrm{Ext}}

\newcommand{\Hom}{\mathrm{Hom}}

\newcommand{\rank}{\mathrm{rank\,}}

\newcommand{\codim}{\mathrm{codim}\,}
\newcommand{\aaa}{\mathbf{a}}
\newcommand{\Red}{\mathcal{R}}
\newcommand{\Euler}[2]{\langle #1\,,\,#2\rangle}
\newcommand{\dv}{\underline{\mathrm{dim}}\,}
\newcommand{\lace}{\mathrm{m}}
\newcommand{\Ind}[1]{\mathrm{Indec}_{#1}}
\newcommand{\Mod}[1]{\text{Mod-$\C#1$}}
\newcommand{\orb}[2]{\overline{#1\cdot#2}}

\newcommand{\Sym}{\mathrm{Sym}}
\newcommand{\Lie}{\mathrm{Lie}}
\newcommand{\del}{\partial}

\newcommand{\AR}[1]{\Gamma_{#1}}

\newcommand{\Schub}{\mathfrak{S}}
\newcommand{\OO}{\mathcal{O}}

\newcommand{\reminder}[1]{\textbf{**#1**}}
\newcommand\junk[1]{}

\begin{document}
\pagestyle{plain}

\title{Kempf collapsing and quiver loci}
\author{Allen Knutson}
\address{Department of Mathematics, University of California, San Diego}
\email{allenk@math.ucsd.edu}

\author{Mark Shimozono}
\address{Mathematics Department\\ Virginia Tech \\ Blacksburg, Virginia}
\email{mshimo@math.vt.edu}

\date{\today}

\maketitle

\begin{abstract}
  Kempf [1976] studied proper, $G$-equivariant maps from
  equivariant vector bundles over flag manifolds to
  $G$-representations $V$, which he called {\em collapsings}.  We give
  a simple formula for the $G$-equivariant cohomology class on $V$, or
  {\em multidegree}, associated to the image of a collapsing: apply a certain
  sequence of divided difference operators to a certain product
  of linear polynomials, then divide by the number of components
  in a general fiber. When that number of components is $1$, we
  construct a desingularization of the image of the collapsing.
  If in addition the image has rational singularities, we can use the
  desingularization to give also a formula for the $G$-equivariant
  $K$-class of the image, whose leading term is the multidegree.

  Our application is to quiver loci and quiver polynomials.
  Let $Q$ be a quiver of finite type (A, D, or E, in arbitrary orientation),
  and assign a vector space to each vertex.
  Let $\Hom$ denote the (linear) space of representations of $Q$ with
  these vector spaces.  This carries an action of $GL$, the product of
  the general linear groups of the individual vector spaces.
  A {\em quiver locus} $\Omega$ is the closure in $\Hom$ of a $GL$-orbit,
  and its multidegree is the corresponding {\em quiver polynomial}.
  Reineke [2004] proved that every ADE quiver locus
  is the image of a birational Kempf collapsing
  (giving a desingularization directly).

  Using Reineke's collapsings, we give formulae for ADE
  quiver polynomials, previously only computed in type A (though in
  this case, our formulae are new).  In the $A$ and $D$ cases quiver
  loci are known to have rational singularities [Bobi\'nski-Zwara
  2002], so we also get formulae for their $K$-classes, which had
  previously only been computed in equioriented type $A$ (and again
  our formulae are new).
\end{abstract}

\setcounter{tocdepth}{1}
{ \tableofcontents}

\section{Introduction}

\renewcommand\Im{\mathop{\mathrm{Im}}}

\subsection{Kempf collapsings}
Let $G$ be a reductive algebraic group, and $P$ a parabolic subgroup.
Let $Y$ be a linear representation of $G$, and $Z\leq Y$ a
$P$-invariant subspace (or more generally, a closed subvariety with at
worst rational singularities). In \cite{K}, Kempf considers
the map
\begin{align*}
  G \times^P Z &\overset{\kappa}{\to} Y & \qquad \hbox{here
    $G\times^P Z :=
    (G\times Z)\big/ \left \{ [g,\vec v] \sim [gp^{-1}, p\vec v] \right\}$}
  \\
  [g, \vec z] &\mapsto g\cdot \vec z &
\end{align*}
which he calls a {\dfn collapsing} of $G\times^P Z$.
This space is the associated fiber bundle with fiber $Z$ over the
homogeneous projective variety $G/P$.  The map $\kappa$ is proper
(it factors as $G\times^P Z \into G/P \times Y \onto Y$),
hence its image $G\cdot Z$ is closed.
When $\kappa$ is birational, $\kappa$ serves as a resolution of
singularities of the variety $G\cdot Z$.

Since $G\cdot Z$ is a $G$-invariant subvariety of $\Hom$, it has a
{\em multidegree} $[G\cdot Z]$, which can be defined as the
associated $G$-equivariant cohomology class (or Chow class) in the
ring $H^*_G(Y)$. Our first result, Theorem \ref{thm:multidegree}
below, is a formula for $[G\cdot Z]$. In Section \ref{sec:mdegs} we
recall the properties we need of multidegrees.

If $T\leq P$ is a maximal torus of $G$, then $H^*_G(Y)$ naturally
includes into $H^*_T(Y) = \Sym^\bullet(T^*)$, so it is enough
to compute $[G\cdot Z]$ as a $T$-equivariant cohomology class.
Since $Z$ is $P$-invariant and hence $T$-invariant,
it too has a $T$-multidegree $[Z] \in \Sym^\bullet(T^*)$.
This $[Z]$ turns out to be particularly simple for $Z$ a linear subspace:
$[Z]$ is the formal product of the $T$-weights in $Y/Z$, each of which lives
in the weight lattice $T^*$ of $T$. If none of the $T$-weights on $Y$ are $0$,
then this product cannot be zero.

\newcommand\Gm{{\mathbb G}_m}

For nonlinear $Z\subseteq Y$, e.g. a union of linear subspaces,
there might still be some cancelation giving $[Z]=0$.  This can't
happen (as follows from Theorem D in \cite{KM}) if all the weights
of $T$ acting on $Y$ live in an open half-space in $T^*$; in this
very common case\footnote{%
  This condition on the weights is not as restrictive as it looks.  If
  $Z\subseteq Y$ is invariant under rescaling (i.e. is the affine cone
  over a projective variety), then we can extend the action of $T$ to
  $T\times \Gm$ where $\Gm$ acts by dilation, and now all the weights
  live in $T^* \times  \{1\}$.  If $Z$ is not already rescaling-invariant,
  we can replace it by the limit subscheme
  $Z' := \lim_{t\to 0} t\cdot Z$, and compute the more refined
  multidegree $[Z'] \in \Sym^\bullet((T\times \Gm)^*)$.
  Afterwards $[Z]$ can be computed as the image of
  $[Z']$ in $\Sym^\bullet(T)^*$ (and this image may indeed be zero).
}
{\em any} closed $T$-invariant scheme $Z\subseteq Y$ has a nonzero
multidegree $[Z]$.

\begin{Theorem}\label{thm:multidegree}
  Let $\kappa : G\times^P Z\to Y$ be a Kempf collapsing,
  where $Lie(P)$ contains all the {\em negative} root spaces.
  Let $d$ be the number of components in a general fiber of $\kappa$.
  Assume that all the weights of $T$ acting on $Y$ live in an open
  half-space in $T^*$.

  Let $m_0 = [Z]$, and construct a sequence of polynomials $m_1, m_2,\ldots$
  by applying divided difference operators
  $\partial_\alpha := \frac{1}{\alpha} (1 - r_\alpha)$ to $m_0$,
  where $\alpha$ varies over the set of simple roots of $G$, and
  $r_\alpha$ acts on $\Sym^\bullet(T^*)$ from the reflection action on $T^*$.
  Don't apply a divided difference operator if the result is $0$,
  and only stop when all $\partial_\alpha$ give the result $0$.

  This process always terminates after the same number of steps
  (namely, $\dim G\times^P Z - \dim G\cdot Z$),
  and the last polynomial in this sequence is $d$ times $[G\cdot Z]$.

  In the case that $\kappa$ is generically finite,
  the sequence of simple roots can be taken to give a reduced
  expression for $w_0 w_0^P$, the product of the long elements of
  the Weyl groups $W$ of $G$ and $W_P$ of $P$ respectively.
  In this case there is an alternate formula
  $$ d\,[G\cdot Z] = \sum_{w\in W^P} w\cdot
  \frac{[Z]}{\prod_{\beta \in \Delta\setminus \Delta_P} \beta} $$
  where $W^P$ is the set of minimal coset representatives in $W/W_P$,
  and $\Delta$ and $\Delta_P$ are the sets of roots of $G$ and $P$
  respectively.
\end{Theorem}

\junk{
  Asking $\kappa$ to be birational is extremely restrictive.
  One sign of the unnaturality of this condition is that even if
  $G\times^P Z \onto G\cdot Z$ is birational,
  the composite $G\times^B Z \onto G\times^P Z \onto G\cdot Z$
  (where $B\leq P$ is a Borel subgroup) will have general fiber $P/B$.
  And given a collapsing $G\times^B Z \onto G\cdot Z$, there may be no
  parabolic $P>B$ preserving $Z$ such that $G\times^P Z \onto G\cdot Z$
  is birational.}

In Section \ref{sec:mdegs} we give an example in which
$\kappa$ does not have connected general fiber.

When $\kappa$ does have connected fibers and $Z$ and $G\cdot Z$ have
rational singularities, we can use the collapsing to compute a more
precise invariant than the multidegree, which is the {\em
$K$-polynomial} $[G\cdot Z]^K_Y$.  Essentially, this is the
numerator of the multigraded Hilbert series of the sheaf
$\OO_{G\cdot Z}$ on $Y$; we recall the precise definition in Section
\ref{sec:mdegs}.

\begin{Theorem}\label{thm:Kpoly}
  Let $\kappa : G\times^P Z\to Y$ be a Kempf collapsing whose general fiber
  is connected (so $d=1$ in the notation of Theorem \ref{thm:multidegree}),
  and assume $Z$ and $G\cdot Z$ have rational singularities.

  Let $m_0 = [Z]^K_Y$, and construct a sequence of Laurent polynomials
  $m_1, m_2,\ldots$ by applying Demazure operators
  $d_\alpha := (1 - \exp(-\alpha))^{-1} (1 - \exp(-\alpha) r_\alpha)$
  to $m_0$, where $\alpha$ varies over the set of simple roots of $G$.
  Don't bother applying any $d_\alpha$ that acts as the identity, and
  only stop when all $d_\alpha$ act as the identity.
  The sequence of simple roots can be taken to give a reduced
  expression for $w_0 w_0^P$.
  This process terminates after finitely many steps
  (namely, $\dim G\times^P Z - \dim G\cdot Z$).
  The last Laurent polynomial in this sequence is the $K$-polynomial
  $[G\cdot Z]^K_Y$. Moreover
  $$ [G\cdot Z]^K_Y = \sum_{w\in W^P} w\cdot
  \frac{[Z]^K_Y}{\prod_{\beta\in \Delta\setminus\Delta_P} (1-\exp(-\beta))}.$$
\end{Theorem}

In the cases where Theorem \ref{thm:Kpoly} applies, it implies
Theorem \ref{thm:multidegree}, by viewing the multidegree as the
lowest-order homogeneous component of the $K$-polynomial.

Kempf worked only with the case that $Z \subseteq Y$ is a linear
subspace (in which case its $K$-polynomial $[Z]^K_Y$ is again a very
simple product over the weights in $Y/Z$), which will also suffice for
our main application. It is frequently the case that the weights of
$T$ on $Y$ are distinct, which implies that there are only finitely many
$P$-invariant linear subspaces $Z$ on which to apply Kempf's construction.

\subsection{Quiver loci}
A {\dfn quiver} $Q=(Q_0,Q_1)$ is a finite directed graph, which
consists of a set $Q_0$ of vertices and a set $Q_1$ of directed
edges or arrows, such that each arrow $a\in Q_1$ has a {\dfn tail}
$ta\in Q_0$ and a {\dfn head} $ha\in Q_0$. A {\dfn representation}
$V$ of $Q$ is a choice of a vector space $V(i)$ for each vertex
$i\in Q_0$ and a linear map $V(a)\in\Hom(V(ta),V(ha))$ for each
arrow $a\in Q_1$. There are obvious notions of isomorphism, direct
sum, and indecomposable, for representations of $Q$. The {\dfn
dimension vector} of $V$ is the map $Q_0\rightarrow \naturals$
defined by $i\mapsto \dim V(i)$. Fix a dimension vector
$d:Q_0\rightarrow\naturals$ and define
$$ GL := GL(Q,d) = \prod_{i\in Q_0} GL(\C^{d(i)}), \qquad
\Hom := \Hom(Q,d) = \prod_{a\in Q_1} \Hom(\C^{d(ta)},\C^{d(ha)}). $$
A typical element of $\Hom$ is denoted $V$, and for $a\in Q_1$ the
$a$ component is denoted $V(a)$. The notation comes from thinking of
$V$ as a functor from the free category on $Q$ to the category {\bf
Vec}. The group $GL$ acts linearly on $\Hom$ by change of basis:
$(g\cdot V)(a) = g(ta) V(a) g(ha)^{-1}$ for all $g\in GL$, $V\in
\Hom$, and $a\in Q_1$.
Two points in $\Hom$ are in the same $GL$-orbit if and only if they
define isomorphic representations of $Q$.
The closures of the $GL$-orbits are called {\dfn quiver loci}.\footnote{%
  The term ``quiver varieties'' is already taken, to refer to the
  hyperk\"ahler quotients $(\Hom \tensor {\mathbb H}) /// GL$.}

\begin{Theorem*}
  \begin{itemize}
  \item \cite{G}
    The action of $GL(Q,d)$ on $\Hom(Q,d)$ has finitely many orbits
    for all dimension vectors $d:Q_0\rightarrow\naturals$, if and only if
    $Q$ is a {\dfn Dynkin quiver}, i.e.
    if the undirected graph underlying $Q$ is a Dynkin diagram of
    type $A_{n\geq 1}$, $D_{n\geq 4}$, $E_6$, $E_7$, or $E_8$.
  \item \cite{LM,BZ}
    For $Q$ of type $A,D$, the quiver loci have rational singularities.
    (To our knowledge the $E$ cases are still open.)
  \item \cite{R}
    If $Q$ is a Dynkin quiver, each quiver locus $\Omega \subseteq \Hom$ is
    the image of a birational linear Kempf collapsing, i.e. there
    exists a parabolic subgroup $P \leq GL$ and a $P$-invariant linear
    subspace $Z\leq \Hom$ such that $GL \times^P Z \onto GL\cdot Z = \Omega$
    is birational.
  \end{itemize}
\end{Theorem*}

In \cite{R}, Reineke constructs each $Z$ explicitly using the {\em
Auslander-Reiten quiver} of $Q$. We recapitulate this construction
precisely in Section \ref{sec:finitequivers}.

Modulo the construction of a certain ordering,
we can state the resulting quiver formulae here.
It is well-known \cite{G} that the indecomposable representations of
a Dynkin quiver $Q$
are in bijection with the set of positive roots $R^+$ of the root
system corresponding to the underlying Dynkin diagram. Fix an ordering
$R^+=\{\beta_1,\beta_2,\dotsc,\beta_N\}$ of the set of positive
roots and write $I_j$ for the (isomorphism class of an)
indecomposable representation of $Q$ corresponding to $\beta_j$.
The correspondence $I_j \leftrightarrow \beta_j$ is determined as follows:
the dimension vector of the indecomposable $I_j$ is given by the expansion
$\beta_j = \sum_{i\in Q_0} \dim I_j(i) \, \alpha_i$
of the corresponding positive root $\beta_j$
in the basis $\{\alpha_i\mid i\in Q_0\}$ of simple roots.

Thus there is a bijection between the $GL$-orbits in
$\Hom=\Hom(Q,d)$ and the direct sums $\bigoplus_{j=1}^N I_j^{\,\oplus m_j}$
where $(m_j\mid j = 1,\ldots,N)$ satisfies the obvious dimension
condition
\begin{align} \label{E:dim}
  \text{for all $i\in Q_0$},\quad
  d(i) = \sum_{j=1}^N d_j(i), \qquad\text{where }
  \ d_j(i) := m_j \dim(I_j(i)).
\end{align}
Fix such a tuple of multiplicities $m = (m_j)$ and let $\Omega_m
\subset \Hom$ be the closure of the corresponding $GL$-orbit.

Based on $m$ we define a parabolic $P_m\subset GL$
and a linear subspace $Z_m\subset\Hom$ as follows.
For each vertex $i\in Q_0$ we divide the sets of row and column indices of
$GL(\C^{d(i)})$ into contiguous subsets of sizes $d_j(i)$ as
$j$ runs from $1$ to $N$.  This defines a standard parabolic subgroup
$P_m\subset GL$ whose $i$th component (for $i\in Q_0$) is the block lower
triangular subgroup of $GL(\C^{d(i)})$ with the given diagonal block sizes.

The decompositions $d(i) = \sum_{j=1}^N d_j(i)$ also induce a block
structure on each component $\Hom(\C^{d(ta)},\C^{d(ha)})$ of $\Hom$,
whose $(j,j')$ block is a $d_j(ta) \times d_{j'}(ha)$ rectangle.
Define the linear subspace $Z_m\subset\Hom$ to be those elements
with zeroes in all blocks strictly above the ``block diagonal". This
$Z_m$ is easily seen to be $P_m$-invariant. Reineke proves that {\em
for certain choices of ordering} (built using reduced words for $w_0$
adapted to the quiver, as spelled out in Section \ref{sec:finitequivers}) 
on $R_+$, the Kempf collapsing
$GL\times^{P_m} Z_m \to \Hom$ is birational to $\Omega_m$.

Let $\{x_k^{(i)}\mid i\in Q_0,\, k \in \{1,\ldots,d(i)\} \}$ be a basis for
the weight lattice $T^*$ of the standard maximal torus $T$
given by the tuples of diagonal matrices in $GL$.
Then the $(k,k')$th matrix entry in the $a$th component of $\Hom$ has
weight $x_k^{(ta)}-x_{k'}^{(ha)}$.
\junk{
  Define the weight of a block in
  $\Hom$ to be the product of the weights of its entries.}

\begin{Theorem} \label{thm:application}
  Let $Q$ be a Dynkin quiver, and $\{\beta_1,\beta_2,\dotsc,\beta_N\}$
  a certain order on the set $R^+$ of positive roots
  (constructed explicitly in Section \ref{sec:finitequivers}).
  Let $\Omega_m$ be a quiver locus, with associated multiplicities $m$,
  parabolic $P_m \leq GL$, and subspace $Z_m \leq \Hom$.

  Then $[\Omega_m]$ may be computed by Theorem \ref{thm:multidegree}
  where $[Z_m]$ is the product of the weights of all blocks in $\Hom$
  that are strictly above the ``block diagonal". For types $A$ and
  $D$ the $K$-polynomial $[\Omega_m]^K_{\Hom}$ may be computed by
  Theorem \ref{thm:Kpoly}, in which $[Z_m]^K_\Hom$ is the product of
  terms of the form $1-e^{-\gamma}$ where $\gamma$ runs over those
  same weights as in $[Z_m]$.
\end{Theorem}

\begin{ex} Let $Q$ be the equioriented $A_n$ quiver:
  
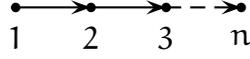
\begin{align*}
    \pspicture(-.5,-.5)(4.5,.5)%
    \psset{arrows=->,arrowsize=3pt 3}%
    \psdots(0,0)(1,0)(2,0)(3,0)%
    \psline[arrows=->](0,0)(1,0)%
    \psline[arrows=->](1,0)(2,0)%
    \psline[arrows=->,linestyle=dashed](2,0)(3,0)%
    \rput(0,-.4){$1$}%
    \rput(1,-.4){$2$}%
    \rput(2,-.4){$3$}%
    \rput(3,-.4){$n$}%
    \endpspicture
  \end{align*}
  The simple roots of $A_n$ are given by
  $\alpha_{ij}=\alpha_i+\alpha_{i+1}+\dotsm+\alpha_j$ for $1\le i\le
  j\le n$. Write $I_{ij}$ for the indecomposable representation of
  $A_n$ corresponding to $\alpha_{ij}$. A suitable ordering on the
  indecomposables is given by
  $I_{11},I_{12},I_{22},I_{13},I_{23},I_{33},\dotsc,I_{nn}$. Let us
  consider the specific example $n=3$, $d=(2,3,2)$, and $\Omega$ given
  by the $GL$-orbit closure of $I_{12}^{\oplus 2} \oplus I_{23}\oplus I_{33}$.
  Geometrically, $\Omega$ is defined by requiring the map
  $V(2) \to V(3)$ to have rank $\leq 1$, and the composite map $V(1)\to V(3)$
  to vanish.

  The decompositions $d(i) = \sum_{j=1}^N d_j(i)$ are
  $$ d(1) = 2+0+0, \quad d(2) = 2+1+0, \quad d(3) = 0+1+1, $$
  so the parabolic $P_m\subset GL$ and the linear subspace
  $Z_m\subset\Hom$ take the following form:
  $$
    P_m = \left\{\left(\begin{pmatrix} *&* \\ *&*
        \end{pmatrix},\begin{pmatrix}
          *&*&0 \\ *&*&0 \\ *&*&*
        \end{pmatrix},
        \begin{pmatrix} *&0 \\ *&*
        \end{pmatrix}\right)\right\}, \qquad
    Z_m = \left\{\left(
        \begin{pmatrix}
          *&*&0\\ *&*&0
        \end{pmatrix}
        ,
        \begin{pmatrix}
          0&0 \\ 0&0 \\ *&0
        \end{pmatrix}
      \right)\right\}.
  $$
  \junk{
  \begin{align*}
    P_m &= \left\{\left(\begin{pmatrix} *&* \\ *&*
        \end{pmatrix},\begin{pmatrix}
          *&*&0 \\ *&*&0 \\ *&*&*
        \end{pmatrix},
        \begin{pmatrix} *&0 \\ *&*
        \end{pmatrix}\right)\right\} \\
    Z_m &= \left\{\left(
        \begin{pmatrix}
          *&*&0\\ *&*&0
        \end{pmatrix}
        ,
        \begin{pmatrix}
          0&0 \\ 0&0 \\ *&0
        \end{pmatrix}
      \right)\right\}
  \end{align*}
  To see this plainly, for the $1\in Q_0$ component of $GL$ the
  diagonal blocks have sizes $2,0,0$, for the $2\in Q_0$ component the
  blocks have sizes $2,1,0$, and for the $3\in Q_0$ component the
  blocks have sizes $0,1,1$. This yields the form of $P_m$ depicted
  above.

  For the $(1,2)\in Q_1$ component of $\Hom$, the above list of
  indecomposables yield ``diagonal" blocks of shapes $2\times 2,
  0\times 1,0\times 0$. Putting submatrices of these shapes corner to
  corner and putting arbitrary entries on or below the ``block
  diagonal" one obtains the $(1,2)$-component of $Z_m$ as depicted
  above. For the $(2,3)\in Q_1$ component the diagonal blocks have
  shapes $2\times0,1\times1,0\times 1$.
  }

  By Theorem \ref{thm:multidegree}
  \begin{align*}
    [Z_m] &= (x^{(1)}_1-x^{(2)}_3)(x^{(1)}_2-x^{(2)}_3)\
    (x^{(2)}_1-x^{(3)}_1)(x^{(2)}_1-x^{(3)}_2)
    (x^{(2)}_2-x^{(3)}_1)(x^{(2)}_2-x^{(3)}_2) (x^{(2)}_3-x^{(3)}_2) \\
    [\Omega] &=
    \partial_{x^{(2)}_1-x^{(2)}_2} \,\partial_{x^{(2)}_2-x^{(2)}_3}
    \,\partial_{x^{(3)}_1-x^{(3)}_2}   \, [Z_m]
  \end{align*}
  By Theorem \ref{thm:Kpoly}
  \begin{align*}
    [Z_m]^K_\Hom &=
    (1-e^{-x^{(1)}_1+x^{(2)}_3})(1-e^{-x^{(1)}_2+x^{(2)}_3})
    (1-e^{-x^{(2)}_1+x^{(3)}_1})(1-e^{-x^{(2)}_1+x^{(3)}_2}) \\
    &\,\,\,\,\,\,\,(1-e^{-x^{(2)}_2+x^{(3)}_1})(1-e^{-x^{(2)}_2+x^{(3)}_2})
    (1-e^{-x^{(2)}_3+x^{(3)}_2}) \\
    [\Omega]^K_\Hom &=
    d_{x^{(2)}_1-x^{(2)}_2}\, d_{x^{(2)}_2-x^{(2)}_3}\, 
    d_{x^{(3)}_1-x^{(3)}_2}\,     [Z_m]^K_\Hom.
  \end{align*}
We shall work out the multidegree calculation explicitly. Let
$a_i=x^{(1)}_i$, $b_i=x^{(2)}_i$, and $c_i=x^{(3)}_i$. We shall use
the following properties of $\partial_\alpha$:
$\partial_\alpha(f)=0$ if $r_\alpha(f)=f$, and
$\partial_\alpha(fg)=\partial_\alpha(f)g+r_\alpha(f)
\partial_\alpha(g)$. In particular if $r_\alpha(f)=f$ then
$\partial_\alpha(fg)=f\partial_\alpha(g)$.

Using the notation $\partial^a_i = \partial_{a_i-a_{i+1}}$ (and
similarly for $b,c$) we have
\begin{align*}
  [\Omega] &= \partial^b_1 \partial^b_2 \partial^c_1
  (a_1-b_3)(a_2-b_3)(b_1-c_1)(b_1-c_2)(b_2-c_1)(b_2-c_2)(b_3-c_2) \\
  &= \partial^b_1 \partial^b_2
  (a_1-b_3)(a_2-b_3)(b_1-c_1)(b_1-c_2)(b_2-c_1)(b_2-c_2) \\
&= \partial^b_1 [ (a_2-b_3)(b_1-c_1)(b_1-c_2)(b_2-c_1)(b_2-c_2) \\
&+ (a_1-b_2)(b_1-c_1)(b_1-c_2)(b_2-c_1)(b_2-c_2) \\
&+ (a_1-b_2)(a_2-b_2)(b_1-c_1)(b_1-c_2)(b_2-c_2) \\
&+ (a_1-b_2)(a_2-b_2)(b_1-c_1)(b_1-c_2)(b_3-c_1) ] \\
&= [0 ]+ [(b_1-c_1)(b_1-c_2)(b_2-c_1)(b_2-c_2)] \\
&+ [(a_2-b_2)(b_1-c_1)(b_1-c_2)(b_2-c_2) \\
&+ (a_1-b_1)(b_1-c_1)(b_1-c_2)(b_2-c_2) \\
&+ (a_1-b_1)(a_2-b_1)(b_1-c_2)(b_2-c_2)] \\
&+ [ (a_2-b_2)(b_1-c_1)(b_1-c_2)(b_3-c_1) \\
&+ (a_1-b_1)(b_1-c_1)(b_1-c_2)(b_3-c_1) \\
&+ (a_1-b_1)(a_2-b_1)(b_1-c_2)(b_3-c_1) \\
&+ (a_1-b_1)(a_2-b_1)(b_2-c_1)(b_3-c_1) ].
\end{align*}
We check this against the component formula \cite[Cor. 6.17]{KMS},
which is a sum over three minimal length lacing diagrams
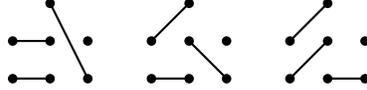
\begin{align*}
  \psset{xunit=.5cm,yunit=.5cm}
  \pspicture(0,1)(2,3)
  \psdots(0,1)(0,2)(1,1)(1,2)(1,3)(2,1)(2,2)
  \psline(0,1)(1,1)
  \psline(0,2)(1,2)
  \psline(1,3)(2,1)
  \endpspicture \qquad
  \pspicture(0,1)(2,3)
  \psdots(0,1)(0,2)(1,1)(1,2)(1,3)(2,1)(2,2)
  \psline(0,1)(1,1)
  \psline(0,2)(1,3)
  \psline(1,2)(2,1)
  \endpspicture \qquad
  \pspicture(0,1)(2,3)
  \psdots(0,1)(0,2)(1,1)(1,2)(1,3)(2,1)(2,2)
  \psline(0,1)(1,2)
  \psline(0,2)(1,3)
  \psline(1,1)(2,1)
  \endpspicture
\end{align*}
which give the three tuples of partial permutation matrices
\begin{align*}
  \left(\begin{pmatrix}
  1&0&0\\ 0&1&0
  \end{pmatrix},\begin{pmatrix}
    0&0 \\ 0&0 \\ 1&0
  \end{pmatrix}\right)\qquad
  \left(\begin{pmatrix}
  1&0&0\\ 0&0&1
  \end{pmatrix},\begin{pmatrix}
    0&0 \\ 1&0 \\ 0&0
  \end{pmatrix}\right)\qquad
  \left(\begin{pmatrix}
  0&1&0\\ 0&0&1
  \end{pmatrix},\begin{pmatrix}
    1&0 \\ 0&0 \\ 0&0
  \end{pmatrix}\right).
\end{align*}
The formula is then
\begin{align*}
  [\Omega] &=
  \Schub_{123}(a;b)\Schub_{3412}(b;c)+\Schub_{132}(a;b)\Schub_{3142}(b;c)+\Schub_{231}(a;b)\Schub_{1342}(b;c)
  \\
  &= [(b_1-c_1)(b_1-c_2)(b_2-c_1)(b_2-c_2)] \\
  &+ [(a_1+a_2-b_1-b_2)((b_1-c_1)(b_1-c_2)(b_2+b_3-c_1-c_2))] \\
  &+ [(a_1-b_1)(a_2-b_1))((b_2-c_1)(b_3-c_1)+(b_1-c_2)(b_3-c_1)+(b_1-c_2)(b_2-c_2))]
\end{align*}
using, say, the pipe dream formula \cite{FK} \cite[Thm. 5.3]{KMS} to
evaluate the double Schubert polynomials $\Schub_w(x;y)$.
\end{ex}


The multidegrees of quiver loci are particularly important for studying the
singularities of composites of differential mappings (see \cite{BF,FR,BFR}
and the references therein).

Until now, the only formulae for these multidegrees were in type $A$.
The first such formula was in \cite{BF}, and applied only to the case
that the directed arrows are all oriented the same direction.
This type $A$ formula has been improved in three ways:
it has been made manifestly positive in an appropriate sense,
the $K$-polynomial has been computed \cite{KMS},
and the orientation has been generalized \cite{BR}.
Some of these have been combined: the $K$-polynomial has been computed
positively \cite{Bu,Mi}, and the multidegree has been computed
positively for arbitrary orientations \cite{BR}.

Using Theorems \ref{thm:multidegree} and \ref{thm:Kpoly}, and the
rationality of the singularities (from \cite{BZ}), we give the first
formulae for
\begin{itemize}
\item the multidegrees of type $D$ and $E$ quiver loci,
\item the $K$-polynomials for type $A$ quiver loci
  in non-equioriented cases, and
\item the $K$-polynomials for type $D$ quiver loci.
\end{itemize}
Unfortunately, our formulae are not positive in the senses of
\cite{KMS,Bu,Mi}. Some positivity of the answers is expected on very
general grounds (e.g. Theorem D in \cite{KM}).

\junk{
  A {\dfn quiver locus} is an orbit closure $\Omega =
  \orb{G}{V}\subset\Hom$ for some $V\in \Hom$. It defines a
  $G$-equivariant cohomology class $[\Omega]\in H^*_G(\Hom)$. Let
  $T\subset G$ be the maximal torus of tuples of diagonal matrices
  and $W$ the Weyl group of $G$. There are natural isomorphisms
  $H^*_T(\Hom) \cong H_T^*(\textrm{pt}) \cong \Sym(\Lie(T))$. Fixing
  a basis of $\Lie(T)$ we may identify this ring with the polynomial
  ring $\Z[x_j^i]$ with generators $\{x_j^i\mid i\in Q_0,1\le j\le
  d(i) \}$. $W$ acts by permuting each set of variables $x_j^i$ for
  $i$ fixed, and $H_G^*(\Hom) \cong H_T^*(\Hom)^W$. The
  {\dfn cohomological quiver polynomial} is the $W$-symmetric
  polynomial given by $[\Omega]$.

  The structure sheaf $\OO_{\Omega}$ of the quiver locus
  $\Omega\subset \Hom$ defines an equivariant class
  $[\OO_{\Omega}]\in K_G^*(\Hom)=K_T^*(\Hom)^W$ called the
  {\dfn $K$-theoretic quiver polynomial}. There are natural
  isomorphisms $K_T^*(\Hom)\cong K_T^*(\textrm{pt})\cong R(T)$ where
  $R(T)$ is the representation ring of $T$. In the coordinates
  chosen above, $R(T)$ is the ring of Laurent polynomials generated
  by the elements $e^{\pm x_j^i}$.
}

\section{The Bott-Samelson crank}
\label{sec:geom}

The inductive processes in Theorems \ref{thm:multidegree} and
\ref{thm:Kpoly} have their geometric origin in the Bott-Samelson
crank \cite{BS}. Fix a Borel subgroup $B$ with $P\geq B\geq T$. For
each simple root $\alpha$ of $G$, let $P_\alpha$ be the
corresponding minimal parabolic.
Then if $f: C\to Y$ is a $B$-equivariant map, the space $P_\alpha
\times^B C$ has also a natural $B$-equivariant map to $Y$, which we
will call $P_\alpha \times^B f$. We call this functor (on the
category of $B$-equivariant maps $f:C\to Y$ to a fixed $G$-space)
{\dfn one turn of the Bott-Samelson crank}. By projecting onto the
first factor, we see that the resulting space is a $C$-bundle over
$P_\alpha/B \iso \Pone$, and in particular $\dim P_\alpha\times^B C
= \dim C + 1$. This $C$-bundle is trivial if $f$ is not just $B$-
but $P_\alpha$-equivariant, with the projection onto the $C$ factor
given by the $B$-quotient of the action map $P_\alpha \times C \to
C$; we study this further in Lemma \ref{lem:GeqvtBS} below.

Since we generally turn the crank many times in succession,
using a sequence $\vec \alpha = (\alpha_1,\ldots,\alpha_k)$,
we will denote products of these functors $P_{\alpha_i} \times^B$
by $BS_{\vec\alpha}
:= P_{\alpha_k} \times^B \ldots \times^B P_{\alpha_1} \times^B$.
A space $BS_{\vec\alpha}\cdot pt$ is a {\dfn
  Bott-Samelson manifold}. The natural $G$-space for a point to map to
$B$-equivariantly is $G/B$, so each Bott-Samelson manifold comes with a
{\dfn Bott-Samelson map} to $G/B$.

Seeing a Bott-Samelson manifold as a free quotient by $B$ on the right of
$P_{\alpha_k} \times^B \ldots \times^B P_{\alpha_1}$, any Bott-Samelson
manifold tautologically carries a principal $B$-bundle. It is sometimes
useful to see the space $BS_{\vec\alpha} Z$ as the associated
$Z$-bundle over the Bott-Samelson manifold $BS_{\vec\alpha}\cdot pt$.

\begin{Lemma}\label{lem:GeqvtBS}
  Let $G$ act on two varieties $C,Y$ (which need not be linear),
  and let $f: C\to Y$ be a $G$-equivariant map. Let $\alpha_1,\ldots,\alpha_j$
  be a sequence of simple roots.

  Then the general fibers of $BS_{\vec\alpha} f : BS_{\vec\alpha} C\to Y$
  have the same number of components as the general fibers of $f$.
\end{Lemma}

\begin{proof}
  Consider the diagram
  $$
  \begin{array}{cccc}
    (BS_{\vec\alpha} \cdot pt) \times C
    &\widetilde\To &  BS_{\vec\alpha} C&  \\
    \downarrow  & & \downarrow& \\
    C & \stackrel{f}{\To} & Y&
  \end{array}
  $$
  The left vertical map is projection onto the second factor,
  and the right vertical map is $BS_{\vec\alpha} f$. If the top map is
  $([p_k,\ldots,p_1],c) \mapsto [p_k,\ldots,p_1,p_1^{-1}\cdots p_k^{-1} c]$,
  which is easily seen to be well-defined and an isomorphism, then
  the diagram commutes.

  We can now study the right-hand map $BS_{\vec\alpha} f$
  by reversing the isomorphism on the top of the diagram. The fibers
  of the map from the northwest corner to the southeast are just
  products of the fibers of $f$ with Bott-Samelson manifolds, which
  are connected.
\end{proof}

\begin{Proposition}\label{prop:BSshrinkdim}
  Let $G$ act on a scheme $Y$, and
  $\iota:Z\into Y$ be the inclusion of a $B$-invariant subvariety.
  (In fact we may as well replace $Y$ by the subvariety $G\cdot Z$.)
  Let $\mu : G\times^B Z \to G\cdot Z$ be
  the projective map $[g,z] \mapsto g\cdot z$.

  Then there exists a sequence of simple roots $(\alpha_1,\ldots,\alpha_k)$,
  such that $BS_{\vec\alpha}\iota$ is surjective and generically finite,
  and its degree is the number of
  components in a general fiber of the map $\mu$.
\end{Proposition}

\junk{
The map $BS_{\vec\alpha}\iota$ will serve later as a 
generically-finite-to-one replacement for $\mu$.
}

\begin{proof}
  \junk{
    Following \cite{K}, we observe that $\mu$
    factors as $G\times^B Z \into G/B \times G\cdot Z \onto G\cdot Z$
    and hence is projective.}

  We will show there exist {\em two} sequences of simple roots
  $(\alpha_1,\ldots,\alpha_k)$, $(\beta_1,\ldots,\beta_j)$ and a
  natural commutative diagram
  $$
  \begin{array}{ccccc}
    BS_{\vec\alpha} Z &\stackrel{BS_{\vec\alpha} \iota} \To&
    G\cdot Z   &\stackrel{\mu}\From& G\times^B Z \\
    &&\uparrow&&\uparrow \\
    && BS_{\vec\beta} (G\cdot Z) &\From& BS_{\vec\beta} BS_{\vec\alpha} Z
  \end{array}
  $$
  in which all maps are onto,
  the map $BS_{\vec\beta} BS_{\vec\alpha} Z \to G\times^B Z$
  is generically $1$:$1$, and
  the map $BS_{\vec\alpha}\iota$ is generically finite to one.
  {}From this diagram we will derive the
  conclusions of the proposition.

  Let $Z_0 = Z$.
  Since the subgroups $\{P_\alpha\}$ generate $G$,
  for each $i$ either $Z_i$ is $G$-invariant or we may pick a simple root
  $\alpha_i$ such that $Z_i$ is not $P_{\alpha_i}$-invariant.
  Define $Z_i := P_{\alpha_i} \cdot Z_{i-1}$.

  Since $Z_i$ is the image of the proper map
  $P_{\alpha_i}\times^B Z_{i-1} \onto Z_i$,
  and by inductive assumption $Z_{i-1}$ is closed, reduced, and irreducible,
  we find $Z_i$ is too. Since $Z_{i-1}$ was not $P_{\alpha_i}$-invariant,
  $Z_i \supset Z_{i-1}$ and $\dim Z_i = \dim Z_{i-1} + 1 = \dim Z + i$.
  Obviously $G\cdot Z_i = G\cdot Z$, so $Z_i \subseteq G\cdot Z$ is
  only $G$-invariant if $Z_i = G\cdot Z$. Hence this process stops when
  $\dim G\cdot Z = \dim Z_k = \dim Z + k$, i.e. $k = \dim G\cdot Z-\dim Z$.

  The map $BS_{\vec\alpha}\iota: BS_{\vec\alpha} Z \to Z_k = G\cdot Z$
  is surjective. By dimension count it is generically finite-to-one.

  \junk{
    Since the map $(\prod_{i=1}^k P_{\alpha_i} \times^B) Z \onto Z_k$
    factors through $(\prod_{i=1}^k P_{\alpha_i} \times^B) Z \onto
    (P_{\alpha_k}\cdots P_{\alpha_1}) \times^B Z$, this latter map is
    also finite-to-one. As explained in \cite{blahblahblah}
    this forces the list of reflections
    $(r_{\alpha_i})$ to be a reduced word in the Weyl group of $G$.
  }

  To construct the sequence $(\beta_j)$, consider the $B$-equivariant map
  $\{pt\} \to G/B$ taking a point to the identity coset, and apply
  $BS_{\vec\alpha}$ to that. The result $BS_{\vec\alpha} \cdot pt \to G/B$
  is a {\dfn Demazure-Hansen resolution} \cite{De,Ha} of a Schubert variety
  in $G/B$, where the source is a Bott-Samelson manifold. Now select
  $(\beta_i)$ following the same procedure as was used above, to construct
  a finite-to-one map $BS_{\vec\beta} BS_{\vec\alpha}\cdot pt \to G/B$.
  In fact the resulting map is generically $1$:$1$ \cite{BS}.
  This obviously extends to a map of $B$-bundles,
  and our map $BS_{\vec\beta} BS_{\vec\alpha} Z \to G \times^B Z$
  is the corresponding map of associated $Z$-bundles.
  Consequently it too is generically $1$:$1$.

  To finish setting up the diagram, define
  $$
  \begin{array}{rclrcl}
    \mu : G\times^B Z &\onto& G\cdot Z
    &\qquad [g,z] &\mapsto& g\cdot z \\
    BS_{\vec\beta} BS_{\vec\alpha} Z &\onto& G\times^B Z
    &\qquad [g_1,\ldots,g_{j+k},z] &\mapsto& [g_1 g_2\cdots g_{k+j},z] \\
    BS_{\vec\beta} BS_{\vec\alpha} Z &\onto& BS_{\vec\beta} Z_k
    =  BS_{\vec\beta} (G\cdot Z)
    &\qquad [g_1,\ldots,g_{j+k},z]
    &\mapsto& [g_1 g_2\cdots g_j,g_{j+1}\cdots g_{j+k} \cdot z]
  \end{array}
  $$
  which are all visibly onto and define the commuting square above. It remains
  to prove our claims about these maps.

  Applying Lemma \ref{lem:GeqvtBS} to the maps
  $BS_{\vec\alpha} \iota : BS_{\vec\alpha} Z\to G\cdot Z$,
  $BS_{\vec\beta} BS_{\vec\alpha} \iota :
  BS_{\vec\beta} BS_{\vec\alpha} Z\to G\cdot Z$,
  we see that the general fiber of
  $BS_{\vec\beta} BS_{\vec\alpha} \iota :
  BS_{\vec\beta} BS_{\vec\alpha} Z\to G\cdot Z$
  has the same number of connected components as the general fiber of
  $BS_{\vec\alpha} \iota : BS_{\vec\alpha} Z\to G\cdot Z$,
  which (since it is generically finite-to-one) is just its degree.
\end{proof}

In the case $Z=pt$,
the following is a standard result about Bott-Samelson manifolds for
partial flag manifolds. 

\begin{Lemma}\label{lem:BSforGmodP}
  Let $Z$ be a $B$-space, and $\vec \alpha$ a list of simple roots whose
  corresponding reflections $(r_{\alpha_i})$ give a reduced word for the
  Weyl group element $w_0 w_0^P$ where $w_0$ is the long element of
  $G$'s Weyl group and $w_0^P$ the long element of $P$'s. Then the
  map $BS_{\vec\alpha} \cdot Z \to G\times^P Z$ (constructed by applying
  $BS_{\vec\alpha}$ to the inclusion $Z \iso P\times^P Z \into G\times^P Z$
  of the fiber over the basepoint) is a birational isomorphism.
\end{Lemma}

\begin{proof}
  These two spaces are $Z$-bundles, and the map takes fibers to fibers;
  as such it is equivalent to check that $BS_{\vec\alpha}\cdot pt\to G/P$
  is a birational isomorphism. Writing this as a composite
  $$ BS_{\vec\alpha}\cdot pt\to G/B \onto G/P, $$
  the first map is birational, by the assumption of reducedness,
  to the (opposite) Schubert variety $\overline{B w_0 w_0^P B}/B$.
  The fiber over $gP \in G/P$ of the second map is
  $gP/B = g\overline{B w_0^P B}/B$.
  Hence the fiber over $gP$ of the composite is the intersection
  $$ g\overline{B w_0^P B}/B \cap \overline{B w_0 w_0^P B}/B $$
  which for generic $g$ is a point, since the $w_0$ makes these
  opposed Schubert varieties.
\end{proof}

\section{Multidegrees, $K$-polynomials, and the proofs of Theorems
  \ref{thm:multidegree} and \ref{thm:Kpoly}} \label{sec:mdegs}

\subsection{Multidegrees and the proof of Theorem \ref{thm:multidegree}.}

Let a torus $T$ act on a vector space $Y$. To each $T$-invariant
subscheme $Z\subseteq Y$, we can associate a {\dfn multidegree} $[Z]_Y$
living in the symmetric algebra on the weight lattice $T^*$ of $T$,
satisfying the following properties:
\begin{enumerate}
\item If $Z=Y=\{0\}$, then $[Z]_Y=1$.
\item If as a cycle $Z = \sum_i m_i Z_i$, where the $\{Z_i\}$ are
  varieties occurring with multiplicities $\{m_i\}$, then
  $[Z]_Y = \sum_i m_i [Z_i]_Y$.
\item If $H \leq Y$ is a $T$-invariant hyperplane, and $Z$ is a variety, then
  \begin{enumerate}
  \item if $Z \not\subseteq H$, then $[Z]_Y = [Z\cap H]_H$, but
  \item if $Z \subseteq H$, then $[Z]_Y = [Z]_H \cdot wt(Y/H)$, where
    $wt(Y/H) \in T^*$ is the $T$-weight on the line $Y/H$.
  \end{enumerate}
\end{enumerate}
The multidegree generalizes the notion of degree of a projective variety
$\PP Z \subseteq \PP Y$. If $T$ is just a circle acting on $Y$ by rescaling,
and $Z$ is the affine cone (hence $T$-invariant) over a projective variety
$\PP Z$, then $[Z]_Y = (\deg \PP Z) a^{\codim_Y Z}$ where $a$ is the
generator of $T^*$. Multidegrees (in $\Sym(T^*)$) are a special case
of equivariant Chow classes (in $A_T(Y)$); since $Y$ is equivariantly
contractible we have $A_T(Y) \iso A_T(pt) \iso \Sym(T^*)$.

It is easy to see that properties (1)-(3) characterize multidegrees.
One can show existence in several ways, one being through multigraded
Hilbert series, as in the next section. Multidegrees were introduced by
\cite{Jo}. Our reference for them is \cite{MS}.

We only use three results about them. One that follows immediately
from the properties above is that for $Z\leq Y$ a linear subspace,
$[Z]_Y$ is the product of the weights in $Y/Z$. The second is that
if all the $T$-weights in $Y$ lie in an open half-space, then $[Z]_Y
\neq 0$ for $Z$ nonempty. (This follows from Theorem D in
\cite{KM}, and is also easily derived from the above properties.) 
The third is a technical result in equivariant Chow theory:

\begin{Lemma}\label{lem:Hlocalization}
  Let $Z$ be a $P$-variety and let $A_T(pt)_{frac}$ denote the field of
  fractions of the polynomial ring $A_T(pt)$. Then we have a formula
  in the localization $A_T(G\times^P Z)\tensor_{A_T(pt)} A_T(pt)_{frac}$
  of the equivariant Chow ring $A_T(G\times^P Z)$:
  $$ 1 = \sum_{w\in W^P} w\cdot
  \frac{[Z]_{G\times^P Z}}{\prod_{\beta\in\Delta\setminus\Delta_P} \beta} $$
  where $[Z]_{G\times^P Z} \in A_T(G\times^P Z)$ is the class induced by
  the regularly embedded subvariety $Z$.
\end{Lemma}

\begin{proof}
  As the map $G\times^P Z \onto G/P$ is $T$-equivariant (indeed,
  $G$-equivariant), all the $T$-fixed points in $G\times^P Z$ lie over
  the $T$-fixed points $\{wP : w\in W^P\}$ in $G/P$. So we get inclusions
  $(G\times^P Z)^T \into \bigcup_{w\in W^P} wZ \into G\times^P Z$.

  Then we use the fact, proven in \cite[section 3.2]{Br},
  that the inclusion of fixed points (the composite of the two above)
  induces an injective pullback $A_T(G\times^P Z) \into A_T((G\times^P Z)^T)$.
  Hence the map $A_T(G\times^P Z) \into
  A_T(\cup_{w\in W^P} wZ) \iso \bigoplus_{w\in W^P} A_T(wZ)$ is injective,
  and to prove the two sides of the formula agree it will suffice to check 
  their images.

  Let $i:Z\into G\times^P Z$ take $z \mapsto [1,z]$.
  Then $i^* [Z]_{G\times^P Z} = i^* i_* 1 = $ the equivariant Euler class
  of the normal bundle of $Z$ inside $G\times^P Z$. This normal bundle
  is the pullback of the normal bundle to the basepoint $P/P \in G/P$,
  hence its equivariant Euler class is the product
  $\prod_{\beta\in\Delta\setminus\Delta_P} \beta$ of the weights in
  the tangent space.

  Applying $i^*$ to both sides of the formula, we therefore get
  $1 = \prod_{\beta\in\Delta\setminus\Delta_P} \beta/
  \prod_{\beta\in\Delta\setminus\Delta_P} \beta$.
  By the Weyl-invariance of both sides, the same confirmation holds
  for the pullback to each $wZ$. Now apply the injectivity above to
  conclude the formula on $G\times^P Z$ itself.
\end{proof}

This has a well-known corollary due to Joseph:

\begin{Corollary}\cite[look in BBM]{Jo}\label{cor:Joseph}
  Let $P_\alpha$ act on $Y$, and $Z$ be a $B$-invariant subscheme.
  Let $d$ be the degree of the map $P_\alpha \times^B Z \to P_\alpha\cdot Z$
  unless $Z$ is $P_\alpha$-invariant, in which case let $d=0$.
  Let $\partial_\alpha$ denote the divided difference operator
  $\frac{1}{\alpha}(1 - r_\alpha)$, acting on $\Sym^\bullet(T^*)$.
  Then
  $$\partial_\alpha [Z]_Y = d\, [P_\alpha \cdot Z]_Y.$$
\end{Corollary}

\begin{proof}
  Let $L$ denote the Levi factor of $P_\alpha$ containing $T$, 
  with semisimple part $L' \iso SL_2$. 
  Then $P_\alpha = L B$, so $P_\alpha \cdot Z = L \cdot Z$.
  Applying Lemma \ref{lem:Hlocalization}, we learn
  $$ \frac{[Z]_{L\times^B Z}}{\alpha}
  + \frac{r_\alpha\cdot [Z]_{L\times^B Z}}{-\alpha} = 1$$
  as elements of $A_T(L\times^{B\cap L} Z)$.
  Let $\kappa$ denote the action map $L\times^{B\cap L} Z \to Y$,
  and apply $\kappa_*$ to both sides:
  $$ \frac{[Z]_Y - r_\alpha\cdot [Z]_Y}{\alpha} = \kappa_*(1). $$
  If $\kappa$ is generically finite of degree $d$, the right-hand side is
  $d\, [L\cdot Z]$, and otherwise $0$.
\end{proof}

(In this corollary we see the reason for $Lie(P)$ to contain all the
{\em negative} root spaces rather than the positive ones; divided difference
operators are usually defined for application to Schubert polynomials,
which come from Schubert varieties that are $B_-$-invariant not $B$-invariant.)

\begin{proof}[Proof of Theorem \ref{thm:multidegree}.]
  Use Proposition \ref{prop:BSshrinkdim} to create a sequence $(\alpha_i)$.
  The condition in Proposition \ref{prop:BSshrinkdim} on $(\alpha_i)$ is that
  $Z_i$ should grow in dimension at each step, which is the condition
  that the $d$ from Corollary \ref{cor:Joseph} is nonzero.
  By the assumption that all the weights of $Y$ lie in a half-space,
  $[P_\alpha \cdot Z]_Y \neq 0$.  Hence the dimension grows if and
  only if $\partial_\alpha$ does not act as zero. So the conditions on
  $(\alpha_i)$ in the theorem's statement match those used in
  Proposition \ref{prop:BSshrinkdim}.

  By Proposition \ref{prop:BSshrinkdim}, the map
  $BS_{\vec\alpha} Z \to Y$ has image $G\cdot Z$.
  The number of components in a general fiber of $G\times^B Z \onto G\cdot Z$
  is the degree of the map $BS_{\vec\alpha} Z \onto G\cdot Z$.
  That degree is in turn the product of the degrees $d_i$ of the maps
  $P_{\alpha_i} \times^B Z_{i-1} \to Z_i$, since
  $BS_{\vec\alpha} Z \onto G\cdot Z$ factors as
  $$ \left(\prod_{i=1}^k P_\alpha \times^B\right) Z \onto
  \left(\prod_{i=1}^{k-1} P_\alpha \times^B\right) Z_1 \onto
  \left(\prod_{i=1}^{k-2} P_\alpha \times^B\right) Z_2 \onto
  \cdots \onto Z_k = G\cdot Z $$
  where the $\{Z_i\}$ are as in the proof of
  Proposition \ref{prop:BSshrinkdim}, and the $j$th map is the associated map
  of bundles over $\left(\prod_{i=1}^j P_\alpha \times^B\right) \cdot pt$
  to the $B$-equivariant map $P_{\alpha_i} \times^B Z_{i-1} \to Z_i$.

  Hence by $k$ applications of Corollary \ref{cor:Joseph},
  $$ d\, [G\cdot Z]_Y = \left(\prod_{i=1}^k d_i\right) [G\cdot Z]_Y
  = \left(\prod_{i=1}^k \partial_{\alpha_i}\right) [Z]_Y. $$

  In the case $\kappa:G\times^P Z \to G\cdot Z$ is generically a finite map,
  we can use Lemma \ref{lem:BSforGmodP} to know that
  for $\vec\alpha$ giving a reduced word for $w_0 w_0^P$, the map
  $BS_{\vec\alpha}\cdot Z \to G\cdot Z$ is also generically finite
  (with the same degree).

  To see the alternate formula, we apply (as in the proof of Joseph's
  Lemma) the pushforward $\kappa_*$ to the equation
  from Lemma \ref{lem:Hlocalization}:
  $$ \kappa_*(1) =
  \sum_{w\in W^P} w\cdot
  \frac{[Z]_Y}{\prod_{\beta\in\Delta\setminus\Delta_P} \beta}. $$
  Since $\kappa$ is generically finite of degree $d$,
  the left-hand side is $\kappa_*(1) = d\, [\Im \kappa]_Y$.
\end{proof}

The first part of this theorem only used Joseph's Lemma (our
Corollary \ref{cor:Joseph}), rather than Lemma
\ref{lem:Hlocalization} directly. This will not be possible in the
proof of Theorem \ref{thm:Kpoly}, where we will use a slightly
different approach.

Kempf assumed a condition on $Z$ that, among other
things, forced the general fiber of a collapsing to be connected.
While his extremely restrictive condition does not hold in our main
application, we will at least have this connectedness,
which is not shared by the following example.

\begin{Example}
Let $G=SL_2(\complexes)$ act on $Y = \lie{sl}_2(\complexes)$ via the
adjoint action, and let $Z = \lie{b}$ be the lower triangular
matrices in $Y$. Let $T$ be the Cartan subgroup of $G$ consisting of
diagonal matrices, and let $P=B$ be the lower triangular matrices in
$G$. Then we run into the problem that the weights of $T$ acting on
$Y$ are $\alpha,0,-\alpha$ where $\alpha$ is the simple root, and do
not all lie in a half-space as required to apply the theorem.

To rescue the example, we enlarge $G$ to $SL_2(\complexes) \times
\complexes^\times$, where the latter circle acts by rescaling on $Y$
and preserves $Z$. Likewise enlarge $T$ and $B$ by this rescaling
circle. Now the weights are $\alpha+a,a,-\alpha+a$ where $a$ is the
generator of the weight lattice of $\complexes^\times$. Recall that
the multidegree $[Z]$ is the product of the weights {\em not}
occurring in $Z$, in this case the one weight $\alpha+a$.

Then the formula gives
$d\, [G\cdot Z] = \partial_\alpha (\alpha+a) = 2$.
And indeed, $G\cdot Z = Y$, so $[G\cdot Z]=1$, while the preimage
in $G\times^B Z$ of a typical diagonal matrix $\diag(t,t^{-1})$ is
$$ \{ (g,z) : \mathrm{Ad}(g)\cdot z = \diag(t,t^{-1}) \} $$
which has $d=2!$ points,
indexed by the permutations of the diagonal entries $t$ and $t^{-1}$.

In the very similar example $G=SL_3(\complexes)$, with $Y,Z,P,B,T$
replaced by their $3\times 3$ counterparts, the general fiber has
$3!$ points. We have
\begin{eqnarray*}
  3!
  &=& \partial_{\alpha_1} \partial_{\alpha_2} \partial_{\alpha_1}
  (a+\alpha_1)(a+\alpha_2)(a+\alpha_1+\alpha_2) \\
  &=& \sum_{w\in S_3} w \cdot \frac
  {(a+\alpha_1)(a+\alpha_2)(a+\alpha_1+\alpha_2)}
  {\alpha_1\, \alpha_2\, (\alpha_1+\alpha_2)}
\end{eqnarray*}
where $\alpha_1,\alpha_2$ are the simple roots of $SL_3(\complexes)$.
\end{Example}

\subsection{$K$-polynomials and the proof of Theorem \ref{thm:Kpoly}.}

A $T$-equivariant coherent sheaf $\mathcal F$ on $Y$ is equivalent to a
$T^*$-graded module $\Gamma$ over $Fun(Y)$. If we assume that the weights
$\{\lambda_i\}$ of $T$ on $Y$ all live in a open half-space of $T^*$,
then each graded piece $\Gamma_\lambda$ is finite-dimensional, and
we can talk about the multigraded Hilbert series $H(\Gamma; t)$.
It is a rational function,
$$ H(\Gamma; t) := \sum_{\lambda \in T^*} \dim(\Gamma_\lambda)\ t^\lambda
= \frac{[\mathcal F]^K_Y}{\prod_{\lambda_i} (1 - t^{\lambda_i})} $$
whose numerator one calls the {\dfn $K$-polynomial} of the sheaf $\mathcal F$.
If $Z$ is a subscheme of $Y$, we will write $[Z]^K_Y$ for the $K$-polynomial
of the structure sheaf of $Z$. It is a function on $T$, i.e. an
element of the Laurent polynomial ring $K_T(Y) \iso K_T(pt)$.

We need some results about $K$-polynomials, corresponding to those we
used about multidegrees. The first, easily calculated from the
Hilbert series definition, is that the $K$-polynomial of a linear
subspace $Z\leq Y$ is the product $\prod (1-t^w)$ where $w$ varies
over the weights of $Y/Z$. The analogue of Lemma
\ref{lem:Hlocalization} is almost word-for-word the same:

\begin{Lemma}\label{lem:Klocalization}
  Let $Z$ be a $P$-variety, and let $K_T(pt)_{frac}$ denote the field of
  fractions of the Laurent polynomial ring $K_T(pt)$. Then we have a formula
  in the localization $K_T(G\times^P Z)\tensor_{K_T(pt)} K_T(pt)_{frac}$
  of the equivariant $K$-ring $K_T(G\times^P Z)$:
  $$ 1 = \sum_{w\in W^P} w\cdot \frac
  {[Z]^K_{G\times^P Z}}{\prod_{\beta\in\Delta\setminus\Delta_P} (1-\exp(-\beta))}$$
  where $[Z]^K_{G\times^P Z} \in A_T(G\times^P Z)$ is the class induced by
  the regularly embedded subvariety $Z$.
\end{Lemma}

\begin{proof}
  Exactly the same proof holds, except that we need localization in
  torus-equivariant algebraic $K$-theory rather than Chow
  \cite[Th\'eor\`eme 2.1]{Th}.
\end{proof}

To apply this formula we need to understand the class
$\kappa_!(1) \in K_T(Y)$. The pushforward $\kappa_!$ in $K$-theory is
defined as the alternating sum of the higher direct images of $\kappa$,
which are difficult to compute in general. 
An especially easy case is when $\kappa$ is
a birational isomorphism, and both spaces have rational singularities; then
$$ \kappa_*(\mathcal O_{G\times^P Z}) = \mathcal O_{G\cdot Z},
\qquad R^i \kappa_*(\mathcal O_{G\times^P Z}) = 0\quad \forall i>0 $$
so $\kappa_!(1) = [G\cdot Z]^K_Y$.

\begin{proof}[Proof of Theorem \ref{thm:Kpoly}.]
  Since $\kappa$ has connected fibers, by Proposition \ref{prop:BSshrinkdim}
  the map $BS_{\vec\alpha}\cdot \iota$ is a birational isomorphism.
  Since $Z$ and $G\cdot Z$ have rational singularities,
  $(BS_{\vec\alpha}\cdot \iota)_!(1) = [G\cdot Z]^K_Y$ as just explained.

  Now we use Lemma \ref{lem:Klocalization} to give a formula for
  $1 \in K_T(BS_{\vec\alpha}\cdot Z)$, and push it forward using
  $(BS_{\vec\alpha}\cdot \iota)_!$, where $\iota:Z\to Y$ is the inclusion.
  Unwinding this formula, we get the first formula claimed.

  (The reason we didn't follow the same induction used in the proof of
  Theorem \ref{thm:multidegree} is that while $Z$ and $G\cdot Z$ have
  rational singularities, we don't know that the intermediate spaces
  constructed in Proposition \ref{prop:BSshrinkdim} do (though this
  seems very likely).)

  The proof of the third formula is exactly the same as in Theorem
  \ref{thm:multidegree}, except that we need to invoke rationality
  of singularities.

  Finally, we prove the second formula from the third, using the map
  $G\times^B Z \onto G\times^P Z$. This is a fibration with fibers $P/B$,
  and the map $\pi:P/B\onto pt$ takes $\pi_!(1) = 1$ (the trivial line
  bundle case of Borel-Weil-Bott). Then we use Lemma \ref{lem:Klocalization}
  to give a formula for $1\in K_T(G\times^B Z)$, which pushes forward to
  the desired formula for $[G\cdot Z]^K_Y$.

  (In $A_T$ rather than $K_T$, the pushforward of $1$ along $P/B\onto pt$
  is zero, which is why there was no analogous formula in Theorem
  \ref{thm:multidegree}.)
\end{proof}


\section{Quiver representations}\label{sec:quiverreps}

A {\dfn representation} $V$ of a quiver $Q$ is a collection
$\{V(i)\mid i\in Q_0\}$ of vector spaces and
$\{V(a)\in\Hom_\C(V(ta),V(ha))\mid a\in Q_1\}$ of linear maps. We
give the reference \cite{GR}.

\subsection{The path algebra $\C Q$} A path of length $m>0$ is a sequence of
arrows $p=a_1a_2\dotsm a_m$ such that $ha_i=ta_{i+1}$ for $1\le i\le m-1$.
The tail and head of the path are given by $tp=ta_1$
and $hp=ha_m$ respectively. One should imagine that one starts at
the vertex $tp=ta_1$ and walks along the arrow $a_1$ to
$ha_1=ta_2$, thence along $a_2$ to $ha_2$, eventually stopping at
$ha_m=hp$. For each $i\in Q_0$ there is a path of length zero also
denoted $i$, with $hi=ti=i$. If $p$ and $p'$ are paths with
$hp=tp'$ then their concatenation $pp'$ is a path. The
{\dfn path algebra} $\C Q$ of the quiver $Q$ is the associative
$\C$-algebra with $\C$-basis given by the set of paths, and
multiplication given by concatenation:
\begin{equation*}
  p \cdot p' = \begin{cases}
  pp' &\text{if $hp=tp'$} \\
  0 &\text{otherwise.}
  \end{cases}
\end{equation*}
$Q_0$ forms a set of orthogonal idempotents for $\C Q$.

\subsection{Modules over $\C Q$} 
Let $\Mod{Q}$ be the category of finite-dimensional
{\em right} $\C Q$-modules. The structure of a module $V\in \Mod{Q}$
is determined as follows. From the action of $Q_0$ there is a direct
sum decomposition $V \cong \bigoplus_{i\in Q_0} V(i)$ where $V(i) :=
V \cdot i$. The map $\dv{V}: Q_0\rightarrow \naturals$ given by
$i\mapsto \dim(V(i))$ is called the {\dfn dimension vector} of $V$.
For $i,j\in Q_0$ and $a\in Q_1$ we have $V \cdot i \cdot a \cdot
j=0$ unless $i=ta$ and $j=ha$. Thus $a$ acts by zero on $V(i)$ for
$i\not=ta$ and defines a linear map $V(a)\in \Hom_\C(V(ta),V(ha))$.
So it is equivalent to work with $\C Q$-modules or with
representations of $Q$.

\begin{rem} \label{rem:rowvec}
  We adopt the convention that matrices act on row vectors.
\end{rem}

\subsection{Quiver loci and quiver polynomials} 
We now change viewpoints, fixing a vector space and the action of the
subalgebra $\C Q_0\subset \C Q$ on it, but letting the rest of the 
$\C Q$-module structure vary.

Fix a dimension vector
$d:Q_0\rightarrow\naturals$. Let
$$\Hom=\Hom(Q,d)=\bigoplus_{a\in Q_1} \Hom_\C(\C^{d(ta)},\C^{d(ha)})$$
be the space of all $\C Q$-module structures on the vector space
$\bigoplus_{i\in Q_0} \C^{d(i)}$ where $\C^{d(i)}$ is the image of
$i\in \C Q_0$. Let $GL=GL(Q,d)=\prod_{i\in Q_0} GL(d(i),\C)$. The
algebraic group $GL$ acts on $\Hom$ by change of basis: $(g\cdot
V)(a) = g(ta) V(a) g(ha)^{-1}$ for all $g\in G$, $V\in \Hom$, and
$a\in Q_1$. It is easy to check that $V,W\in\Hom$ are isomorphic as
elements of $\Mod{Q}$ if and only if they are in the same
$GL$-orbit.

\subsection{Indecomposables and multiplicities}

We want a nice way to index the quiver loci, which are in bijection
with the isomorphism classes in $\Mod{Q}$. Let $\Ind{Q}$ denote the
set of isomorphism classes of indecomposables in $\Mod{Q}$. For
simplicity of notation, we will sometimes write $U$ instead of
$[U]$. For $V\in \Mod{Q}$ and $U\in \Ind{Q}$, define the
multiplicities $\lace_U(V)$ of $V$ by
\begin{equation}
  \label{eq:lace} V \cong \bigoplus_{U\in\Ind{Q}} U^{\oplus \lace_U(V)}.
\end{equation}
The multiplicities $\lace(V)=(\lace_U(V)\mid U\in \Ind{Q})$
determine $V$ up to isomorphism. Let
$\Omega_\lace:=\overline{GL\cdot V}$ for any $V$ with multiplicities
$\lace$. For the equioriented type $A$ quiver the multiplicities
were in \cite{KMS} called the ``lace array''.

\subsection{The Auslander-Reiten quiver}
\label{ss:AR} We recall the definition of the Auslander-Reiten
quiver $\AR{Q}$ associated to the category $\Mod{Q}$ \cite{ARS}.

A map $f$ is {\dfn irreducible} if for all compositions of maps
$f=gh$ with neither $g$ nor $h$ the identity, $g$ is not a split
monomorphism and $h$ is not a split epimorphism.

The {\dfn Auslander-Reiten quiver $\AR{Q}$ of $Q$} is the directed graph
whose vertex set is $\Ind{Q}$ with a directed edge from $[V]$ to
$[W]$ if and only if there is an irreducible map $V\rightarrow W$.

\subsection{Extensions}
\label{SS:ext} For $V,W\in\Mod{Q}$, call $E\in \Mod{Q}$ an {\dfn
extension of $V$ by $W$} if there is a short exact sequence $0
\rightarrow W \rightarrow E\rightarrow V \rightarrow 0$ of $\C
Q$-modules. For each $i\in Q_0$ choose a basis of $E(i) \cong
W(i)\oplus V(i)$ that consists of a basis of $W(i)$ followed by a
basis of $V(i)$ and write the linear maps with respect to this
basis. With our row-vector conventions of Remark \ref{rem:rowvec}, $E(a)$ has
the form
\begin{equation} \label{eq:block}
  E(a) = \begin{pmatrix} W(a) & 0 \\
    * & V(a)
  \end{pmatrix}.
\end{equation}
Let $E(V,W)$ be the set of extensions of $V$ by $W$ with fixed
underlying vector space $V\oplus W$. There is a linear isomorphism
\begin{equation} \label{eq:extmap}
  \bigoplus_{a\in Q_1} \Hom_\C(V(ta),W(ha))\rightarrow E(V,W)
\end{equation}
whose $a$-th component is given by replacing the submatrix $*$ in
\eqref{eq:block} with the element of $\Hom_\C(V(ta),W(ha))$ for $a\in Q_1$.

Say that $E,E'\in E(V,W)$ are equivalent if there is a $\C Q$-module
isomorphism $E\rightarrow E'$ whose restriction to $W$ is the
identity and whose induced map $E/W\rightarrow E'/W$ is the
identity. $\Ext^1_Q(V,W)$ is isomorphic to $E(V,W)$ modulo the above
equivalence (see for example \cite[Thm. 7.21]{Ro}).

\subsection{The canonical resolution} For $V,W\in\Mod{Q}$ let
$\Hom_Q(V,W)$ be the space of right $\C Q$-module homomorphisms from
$V$ to $W$. There is an exact sequence \cite{Ri}
\begin{equation} \label{eq:res}
  \begin{split}
    0 \rightarrow \Hom_Q(V,W) \overset{j}\rightarrow
    \displaystyle{\bigoplus_{i\in Q_0}} \Hom(V(i),W(i))
    &\overset{d_V^W}\rightarrow \displaystyle{\bigoplus_{a\in Q_1}}
    \Hom(V(ta),W(ha)) \\
    &\overset{p}{\rightarrow} \Ext^1_Q(V,W) \rightarrow 0
  \end{split}
\end{equation}
where $j$ is inclusion, $p$ is induced by the map in
\eqref{eq:extmap} and $d_V^W$ is given by
$$(d_V^W(f))_a = V_a f_{ha} - f_{ta} W_a\qquad\text{for $a\in
Q_1$.}$$

The exactness of \eqref{eq:res} gives
\begin{align}
  \label{eq:drank}
  \dim \Hom_\C(V,W) &= \rank d_V^W+\dim \Hom_Q(V,W).
\end{align}

\subsection{The homological form} 
Let $V,W\in \Mod{Q}$. The {\dfn homological form} is defined by
\begin{align*}
  \Euler{V}{W} = \sum_{i\ge0} (-1)^i \dim \Ext_Q^i(V,W).
\end{align*}
The exact sequence \eqref{eq:res} implies that $\Mod{Q}$ is
hereditary (that is, $\Ext_Q^i(V,W)=0$ for $i\ge 2$) and its
exactness gives
\begin{equation} \label{eq:HomExt}
\begin{split}
  \Euler{V}{W} &= \dim \Hom_Q(V,W) - \dim \Ext^1_Q(V,W)  \\
  &= \sum_{i\in Q_0} \dim V(i) \dim W(i) - \sum_{a\in Q_1} \dim
  V(ta) \dim W(ha) \\
  &= \Euler{\dv V}{\dv W}
\end{split}
\end{equation}
where, for dimension vectors $d,d':Q_0\rightarrow \naturals$ we
write
\begin{equation*}
  \Euler{d}{d'} = \sum_{i\in Q_0} d(i)d'(i) - \sum_{a\in Q_1}
  d(ta)d'(ha).
\end{equation*}

\subsection{Codimension and Ext}

By \eqref{eq:HomExt} for $V=W$ and the fact that $\Hom_Q(V,V)$ is
the closure of the stabilizer of $V$ in $GL$, we have
\begin{align*}
  \dim GL - \dim \Hom &= \Euler{V}{V} = \dim \Hom_Q(V,V) - \dim
  \Ext^1_Q(V,V) \\
  &= (\dim GL - \dim GL\cdot V) - \dim \Ext^1_Q(V,V).
\end{align*}
This implies that for $V\in\Hom$, we have
\begin{align*}
  \codim \overline{GL\cdot V} =\dim \Ext^1_Q(V,V).
\end{align*}
Let $\lace$ be a set of multiplicities with $\Omega_\lace \subset
\Hom$. Then
\begin{equation} \label{eq:codim}
  \codim\, \Omega_\lace =
  \sum_{U,W\in \Ind{Q}} \lace_U \lace_W \,\dim \Ext^1_Q(U,W).
\end{equation}

\section{Quivers of finite type}\label{sec:finitequivers}
Let $X_n$ be a simply-laced root system of rank $n$; it is either
$A_n$ for $n\geq 1$, $D_n$ for $n\geq 4$, or $E_n$ for $n=6,7,8$, 
where $n$ is always the number of nodes in the Dynkin diagram. 
We shall also write $X_n$ for the undirected graph given by its Dynkin diagram.
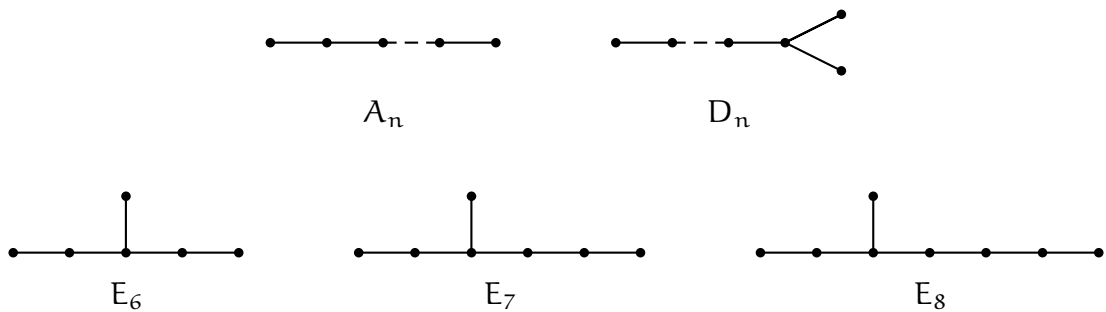
\begin{align*}
\psset{xunit=.75cm,yunit=.75cm}
\begin{matrix}
\pspicture(-.5,-.75)(4.5,.75)%
  \psdots(0,0)(1,0)(2,0)(3,0)(4,0) %
  \psline(0,0)(2,0) %
  \psline(3,0)(4,0) %
  \psline[linestyle=dashed](2,0)(3,0)%
\endpspicture \\
A_n
\end{matrix} \qquad
\begin{matrix}
\pspicture(-.5,-.75)(4.5,.75)%
  \psdots(0,0)(1,0)(2,0)(3,0)(4,.5)(4,-.5) %
  \psline(0,0)(1,0)%
  \psline[linestyle=dashed](1,0)(2,0)%
  \psline(2,0)(3,0)(4,.5)(3,0)(4,-.5)%
\endpspicture \\
D_n
\end{matrix}
\end{align*}
\begin{align*}
\psset{xunit=.75cm,yunit=.75cm}
\begin{matrix}
\pspicture(-.5,-.25)(4.5,1.25)%
  \psdots(0,0)(1,0)(2,0)(3,0)(4,0)(2,1) %
  \psline(0,0)(4,0)
  \psline(2,0)(2,1)%
\endpspicture \\
E_6
\end{matrix}\qquad
\begin{matrix}
\pspicture(-.5,-.25)(5.5,1.25)%
  \psdots(0,0)(1,0)(2,0)(3,0)(4,0)(5,0)(2,1) %
  \psline(0,0)(5,0) %
  \psline(2,0)(2,1)%
\endpspicture \\
E_7
\end{matrix}\qquad
\begin{matrix}
\pspicture(-.5,-.25)(6.5,1.25)%
  \psdots(0,0)(1,0)(2,0)(3,0)(4,0)(5,0)(6,0)(2,1) %
  \psline(0,0)(6,0) %
  \psline(2,0)(2,1)%
\endpspicture \\
E_8
\end{matrix}
\end{align*}

An {\dfn orientation} of an undirected multigraph is a quiver
obtained by choosing directions for the edges of the undirected
graph. Orientations of the Dynkin diagrams of simply-laced root
systems are called {\dfn Dynkin} quivers.

A quiver $Q$ is of {\dfn finite type} if, for every dimension vector
$d:Q_0\rightarrow\naturals$, there are finitely many isomorphism
classes of representations of $Q$ with dimension vector $d$. By
Gabriel's Theorem \cite{G} a quiver is of finite type if and only if
it is Dynkin. In this section we shall assume that $Q$ is Dynkin.

\subsection{Dimension vectors and roots} We recall some
well-known results of Gabriel. There is a bijection from $Q_0$ to
the set of simple roots of $X_n$ given by $i\mapsto \alpha_i$. Any
dimension vector $d:Q_0\rightarrow\naturals$ may be viewed as an
element of the positive cone of roots $\bigoplus_{i\in Q_0}
\naturals \alpha_i$, namely, $\sum_{i\in Q_0} d(i) \alpha_i$. Let
$R^+$ be the set of positive roots of $X_n$.\footnote{We use this
notation to distinguish the root system of $X_n$ with that of the
group $GL$.} There is a bijection $\Ind{Q}\rightarrow R^+$ given by
$U\mapsto \dv U$. $U$ is indecomposable if and only if $\Euler{\dv
U}{\dv U}=1$.

\subsection{Dynkin quivers and orders on $R^+$}
Let $s_i$ denote a simple reflection for the Weyl group $W(X_n)$ of $X_n$%
\footnote{Again this notation is to distinguish $s_i$ from the
  reflection $r_i$ in the Weyl group of $GL$.}
and let $w_0\in W(X_n)$ be the longest element. For $w\in
W(X_n)$ let $\Red(w)\subset Q_0^{\ell(w)}$ denote the set of reduced
words for $w$.

Given an orientation $Q$ of $X_n$ and a vertex $i\in Q_0$, let $s_i
Q$ be the orientation of $X_n$ given by reversing all arrows with
head $i$. Say that a reduced word $\aaa=a_1a_2\dotsm\in\Red(w_0)$ is
{\dfn adapted} to the orientation $Q$ of $X_n$ if $a_j$ is a sink 
(the tail of no arrow) 
in $s_{a_{j-1}}\dotsm s_{a_2}s_{a_1} Q$ for all $j$. By \cite{BGP}, for
every orientation $Q$ of $X_n$, there is a reduced word
$\aaa\in\Red(w_0)$ that is adapted to $Q$.

Each reduced word $\aaa=a_1a_2\dotsm\in\Red(w_0)$ defines a linear
ordering on $R^+$ given by
\begin{equation} \label{eq:rootorder}
  \gamma_1<\gamma_2<\dotsm
\end{equation}
where
\begin{equation} \label{eq:rootdef}
  \gamma_j=s_{a_1} \dotsm  s_{a_{j-1}} (\alpha_{a_j}).
\end{equation}

\subsection{Auslander-Reiten quiver reprise}
There is a combinatorial recipe for the Auslander-Reiten quiver
$\AR{Q}$ of a quiver $Q$ that is an orientation of a Dynkin diagram
$X_n$ of type ADE. This is well-known to the experts; see \cite{Be,Z}.

The vertices of $\AR{Q}$ shall be drawn in the plane in rows
indexed by the set $Q_0$ and columns indexed by $\Z_{>0}$.

Let $\aaa\in\Red(w_0)$ be adapted to $Q$. Let $\gamma_j\in R^+$ be
defined as in \eqref{eq:rootdef}. 
Let $c_1 = 1$, and $c_j = c_{j-1}$ unless for some $k<j$
with $c_k = c_{j-1}$, $\gamma_k$ is adjacent in $X_n$ to $\gamma_j$;
in this case let $c_j = c_{j-1}+1$. 
The vertex $\gamma_j$ is drawn in row $a_j$ and column $c_j$.
Draw a directed edge from $\gamma_j$ to
$\gamma_k$ if $j<k$, $a_j$ and $a_k$ are adjacent in $X_n$, and $k$
is minimal with this property.

\begin{ex} Let $X_n=D_4$ with orientation $Q$ given below.
  
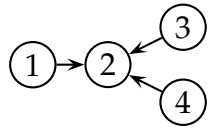
\begin{equation*}
    \pspicture(0,-.75)(2,.75)%
    \cnodeput(2,.5){A1}{3}%
    \cnodeput(0,0){A3}{1}%
    \cnodeput(1,0){A2}{2}%
    \cnodeput(2,-.5){A4}{4}%
    \ncline[arrowscale=1.5]{->}{A1}{A2} %
    \ncline[arrowscale=1.5]{<-}{A2}{A3} %
    \ncline[arrowscale=1.5]{<-}{A2}{A4} %
    \endpspicture
  \end{equation*}
  One reduced word adapted to $Q$ is $213423142341$. The
  corresponding roots have expansions in the simple roots by the
  following matrix.
  \begin{align*}
  \begin{pmatrix}
  \gamma_1 \\ \gamma_2 \\ \gamma_3 \\
\gamma_4 \\ \gamma_5 \\ \gamma_6 \\
\gamma_7 \\ \gamma_8 \\ \gamma_9 \\
\gamma_{10} \\ \gamma_{11} \\ \gamma_{12}  \end{pmatrix} =
\begin{pmatrix}
 0&1&0&0 \\
 1&1&0&0 \\
 0&1&1&0 \\
 0&1&0&1 \\
 1&2&1&1 \\
 1&1&0&1 \\
 0&1&1&1 \\
 1&1&1&0 \\
 1&1&1&1 \\
 0&0&0&1 \\
 1&0&0&0 \\
 0&0&1&0
\end{pmatrix}
\cdot
\begin{pmatrix}
\alpha_1 \\ \alpha_2 \\ \alpha_3 \\ \alpha_4
\end{pmatrix}
\end{align*}
  The Auslander-Reiten quiver $\AR{Q}$ is given by
  
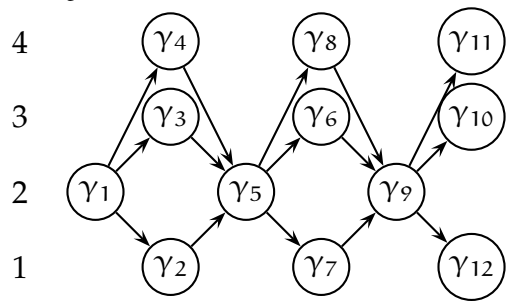
\begin{equation*}
    \pspicture(-1,.75)(5,4.25)%
    \rput(-1,4){4} %
    \rput(-1,3){3} %
    \rput(-1,2){2} %
    \rput(-1,1){1} 
    \cnodeput(0,2){G1}{$\gamma_1$}%
    \cnodeput(1,1){G2}{$\gamma_2$}%
    \cnodeput(1,3){G3}{$\gamma_3$}%
    \cnodeput(1,4){G4}{$\gamma_4$}%
    \cnodeput(2,2){G5}{$\gamma_5$}%
    \cnodeput(3,1){G7}{$\gamma_7$}%
    \cnodeput(3,3){G6}{$\gamma_6$}%
    \cnodeput(3,4){G8}{$\gamma_8$}%
    \cnodeput(4,2){G9}{$\gamma_9$}%
    \cnodeput(5,1){G12}{$\gamma_{12}$}%
    \cnodeput(5,3){G10}{$\gamma_{10}$}%
    \cnodeput(5,4){G11}{$\gamma_{11}$}%
    \ncline[arrowscale=1.5]{->}{G1}{G2} %
    \ncline[arrowscale=1.5]{->}{G1}{G3} %
    \ncline[arrowscale=1.5]{->}{G1}{G4} %
    \ncline[arrowscale=1.5]{->}{G2}{G5} %
    \ncline[arrowscale=1.5]{->}{G3}{G5} %
    \ncline[arrowscale=1.5]{->}{G4}{G5} %
    \ncline[arrowscale=1.5]{->}{G5}{G6} %
    \ncline[arrowscale=1.5]{->}{G5}{G7} %
    \ncline[arrowscale=1.5]{->}{G5}{G8} %
    \ncline[arrowscale=1.5]{->}{G6}{G9} %
    \ncline[arrowscale=1.5]{->}{G7}{G9} %
    \ncline[arrowscale=1.5]{->}{G8}{G9} %
    \ncline[arrowscale=1.5]{->}{G9}{G10} %
    \ncline[arrowscale=1.5]{->}{G9}{G11} %
    \ncline[arrowscale=1.5]{->}{G9}{G12} %
    \endpspicture
  \end{equation*}
  Since nodes $1,3,4$ have no connections in $D_4$, the orders they
  appear in the reduced word $2\,134\,2\,314\,2\,341$ don't affect
  the shape of the Auslander-Reiten quiver.
\end{ex}

\begin{rem} \label{rem:DynkinExt}
  For $Q$ an orientation of the Dynkin diagram of a simply-laced
  root system $X_n$ and $\aaa\in\Red(w_0)$ a reduced word adapted to
  $Q$, let the positive roots (hence the indecomposables) be totally
  ordered as in \eqref{eq:rootorder}. Then for $V,W\in\Ind{Q}$
  we have   \cite{Ri2}
  \begin{align}
  \label{E:Extzero}
    \Ext^1_Q(V,W) &= 0 \qquad&\text{if $V\le W$}.
  \end{align}
\end{rem}

\subsection{The poset of quiver loci in $\Hom$}

\begin{thm} \label{thm:oc} \cite{Bo}
  Let $Q$ be of finite type and $V,W\in\Mod{Q}$
  with $\dv{V}=\dv{W}$. The following are equivalent:
  \begin{enumerate}
  \item $\overline{GL\cdot V} \subset \overline{GL\cdot W}$.
  \item $\dim \Hom_Q(U,V) \le \dim \Hom_Q(U,W)$ for all $U\in\Ind{Q}$.
  \item $\dim \Ext^1_Q(U,V)\ge \dim \Ext^1_Q(U,W)$ for all $U\in\Ind{Q}$.
  \item $\rank d_U^V \ge \rank d_U^W$ for all $U\in\Ind{Q}$.
  \end{enumerate}
\end{thm}
Note that the latter three are equivalent for any quiver $Q$, by
\eqref{eq:HomExt} and \eqref{eq:drank}.

\subsection{The Reineke filtration} We recall a special case of
Reineke's filtration \cite{R}. Let $Q$ be an orientation of a Dynkin
diagram $X_n$ of type ADE, $\aaa\in\Red(w_0)$ adapted to $Q$, with
the associated total order $\le$ on $\Ind{Q}$. We list the elements
of $\Ind{Q}$ in descending order:
$\Ind{Q}=\{\beta_1>\beta_2>\dotsm>\beta_N\}$ where $N=|R^+|$; the
decreasing indexing is for technical convenience related to our
row-vector convention of Remark \ref{rem:rowvec}. For short we write $I_j$ for
the indecomposable instead of $I_{\beta_j}$. Let $V\in\Mod{Q}$,
$d=\dv{V}$, $GL=GL(Q,d)$, $\Hom=\Hom(Q,d)$. Let $V$ have
multiplicities $\lace_j(V)=\lace_{I_j}(V)$ as in \eqref{eq:lace}.
For $1\le j\le N$ write $W_j = I_j^{\oplus \lace_j(V)}$ and $V_j =
W_1 \oplus\dotsm\oplus W_j$. Let $P\subset GL$ be the parabolic
subgroup such that for all $i\in Q_0$, the $i$-th component
$P(i)\subset GL(\C^{d(i)})$ is the stabilizer of $V_j(i)$ for all
$1\le j\le N$. Note that $P$ has Levi factor $L=\prod_{i\in Q_0}
\prod_{j=1}^N GL(W_j(i))\cong \prod_{j=1}^N GL(W_j)$. Let
$Z,Z'\subset Y:=\Hom$ be the coordinate subspaces defined by $Z'(a)
= \bigoplus_{j=1}^N \Hom_{\C}(W_j(ta),W_j(ha))\subset
\Hom_{\C}(V(ta),V(ha))$ and $Z(a)=\bigoplus_{1\le j\le m\le N}
\Hom_{\C}(W_m(ta),W_j(ha))$. For each $a\in Q_1$, $Z'(a)$ is ``block
diagonal" and $Z(a)$ is ``block lower triangular" inside the
matrices\break
$\Hom_{\C}(V(ta),V(ha))$. We claim that
\begin{equation}
  Z = \overline{P \cdot V}
\end{equation}
inside $\Hom$. By \eqref{eq:codim} and \eqref{E:Extzero} we have
$\codim_{\Hom(Q,W_j)} \overline{GL(W_j)\cdot
W_j}=\Ext^1_Q(W_j,W_j)=0,$ or equivalently, $\overline{L\cdot V} =
Z'$. So it suffices to show
\begin{equation}
  Z = \overline{U \cdot Z'}
\end{equation}
where $U$ is the unipotent radical of $P$. But this follows by
induction from the definition of Ext in Subsection \ref{SS:ext}
combined with the fact that by \eqref{E:Extzero} we have
\begin{equation}
\Ext^1_Q(W_p,W_q) = 0\qquad\text{for $p<q$.}
\end{equation}

The linear space $Z$ is the base of our Bott-Samelson induction. Given a
quiver locus $\Omega=\overline{GL\cdot V}\subset \Hom$, we start
with $Z=\overline{P \cdot V}\subset\Hom$. Then $\overline{GL\cdot Z}
=\Omega$. Since $Z$ is a coordinate subspace, $[Z]\in H_T^*(\Hom)$
and $[\OO_Z]\in K_T^*(\Hom)$ have simple product formulae. Applying
Theorem \ref{thm:multidegree} we obtain divided difference formulae
for the multidegree of the quiver locus $\Omega$. By Theorem
\ref{thm:Kpoly}, for quivers of type AD we obtain divided difference
formulae for the $K$-polynomial of $\Omega$.

\begin{Remark} We use an unnecessarily fine filtration. One may
use a directed partition of $R^+$ as defined in \cite{R} to obtain a
coarser filtration of $V$, which leads to a more efficient divided
difference formula.
\end{Remark}

\begin{Example} Let $Q$ be the type $A_2$ quiver with dimension
vector $(m,n)$. Then $\Hom(Q,d)=M_{m\times n}(\complexes)$. For each
$0\le r\le \min(m,n)$ there is a quiver locus $\Omega_r\subset
M_{m\times n}(\complexes)$ given by the determinantal variety of
matrices of rank at most $r$. Using the reduced word $s_2 s_1s_2$ we
have $W_1=I_{\alpha_1}^{\oplus (m-r)}$, $W_2 =
I_{\alpha_1+\alpha_2}^{\oplus r}$, and $W_3 = I_{\alpha_2}^{\oplus
(n-r)}$. The indecomposables can be realized by matrices as follows:
$I_{\alpha_1}$ is a $1\times 0$ matrix, $I_{\alpha_1+\alpha_2}$ can
be taken to be the $1\times 1$ identity matrix, and $I_{\alpha_2}$
is the $0\times 1$ matrix. With respect to bases adapted to the
ordered direct sum $V=W_1\oplus W_2\oplus W_3$, $V\in M_{m\times
n}(\complexes)$ has the $r\times r$ identity matrix in its lower
left corner and zeroes elsewhere. Then $P(1)\subset GL(m)$ and
$P(2)\subset GL(n)$ are block lower triangular with diagonal blocks
of sizes $(m-r,r)$ and $(r,n-r)$ respectively. We have $Z=Z'$; both
are equal to the linear subspace of $M_{m\times n}$ where the bottom
left $r\times r$ submatrix is arbitrary and the other entries are
zero. Let $T(m)\subset GL(m)$ and $T(n)\subset GL(n)$ have weights
$X=(x_1,\dotsc,x_m)$ and $Y=(y_1,\dotsc,y_n)$ respectively. Since
the parabolics $P(1)$ and $P(2)$ are lower triangular, the positive
roots of $P(1)$ have weights $x_j-x_i$ for $1\le i<j\le m$ and those
of $P(2)$ have weight $y_j-y_i$ for $1\le i<j\le n$.

So for $(m,n)=(2,3)$ and $r=1$ we have
\begin{align*}
  [Z] &= (x_1-y_1)(x_1-y_2)(x_1-y_3)(x_2-y_2)(x_2-y_3) \\
  [\Omega] &= \del_{x_1-x_2} \del_{y_1-y_2} \del_{y_2-y_3}  [Z] \\
  &= s_2[X-Y],
\end{align*}
the double Schur polynomial. In general the multidegree is given by
the Giambelli-Thom-Porteous formula $[\Omega_r] = s_{(m-r)\times
(n-r)}[X-Y]$, where the answer is the double Schur polynomial
indexed by the $(m-r)\times (n-r)$ rectangle.
\end{Example}

\section{Beyond $ADE$ quivers}

Let $Q$ be a quiver, $d:Q_0 \to \naturals$ a dimension vector, and
$\Hom$ the associated space of representations. Then as long as $Q$
has no self-loops ($ta=ha$ for some edge $a$), and no repeated edges
($ta=tb$, $ha=hb$ for two edges $a\neq b$), the weights of $T$ on
$\Hom$ are all distinct.

Consequently, there are only finitely many $T$-invariant subspaces
in $\Hom$ (precisely $2^{\dim\Hom}$), and hence {\em only finitely many}
$B$-invariant subspaces $Z$ to which to apply Kempf's construction.
Whereas there may be infinitely many quiver loci.
This (and the fact that quiver loci can have bad singularities
\cite[section 6]{Zw}) suggests that instead of quiver loci, perhaps the
better-behaved objects of study are the $GL$-sweeps of the $B$-invariant
subspaces. From this point of view it is merely an accident (and
Reineke's theorem) that in the $ADE$ case, the two notions coincide.

\end{document}